\documentclass[12pt]{amsart}
\usepackage{amssymb,amsmath,latexsym,amscd,amsfonts,amsrefs,amsthm}
\usepackage{graphics}
\usepackage{epsfig}
\usepackage{psfrag}
\usepackage{cases}
\usepackage{enumerate}
\def\ignore #1 {}

\newtheorem{thm}{Theorem}[section]
\newtheorem{prop}[thm]{Proposition}
\newtheorem{lem}[thm]{Lemma}
\newtheorem{cor}[thm]{Corollary}

\theoremstyle{definition}
\newtheorem{dfn}[thm]{Definition}
\newtheorem{rem}[thm]{Remark}
\newtheorem{example}[thm]{Example}

\def\hpic #1 #2 {\mbox{$\begin{array}[c]{l} \epsfig{file=#1,height=#2} \end{arr\
ay}$}}
\def\vpic #1 #2 {\mbox{$\begin{array}[c]{l} \epsfig{file=#1,width=#2} \end{arra\
y}$}}

\def\B{{\rm B}}

\def\N{\mathbb N}

\def\R{\mathbb R}
\def\P{\mathcal P}
\def\O{\mathcal O}
\def\Q{\mathbb Q}
\def\S{\mathcal S}

\def\OP{{\rm OP}}

\def\A{\mathcal A}

\def\op{{\rm op}}

\def\uex{{\rm uex}}
\def\lex{{\rm lex}}
\def\ex{{\rm ex}}
\def\Dar{{\rm Dar}}
\def\BDar{{\rm BDar}}
\def\Aut{{\rm Aut}}

\def\ev{{\rm ev}}

\def\ev{{\rm ev}}
\def\2{{\mathbf 2}}
\def\1{{\mathbf 1}}

\def\id{{\rm id}}
\def\Dom{{\rm dom}}
\def\im{{\rm im}}
\def\Par{{\rm Par}}
\def\Int{{\rm Int}}
\def\int{{\rm int}}
\def\hto{\rightharpoonup}

\def\proof{\noindent\emph{Proof:\,\,\,}}

\begin{document}

\title{Darboux Calculus}

\author{Marco Aldi and Alexander McCleary}

\begin{abstract} 
We introduce a formalism to analyze partially defined functions between ordered sets. We show that our construction provides a uniform and conceptual approach to all the main definitions encountered in elementary real analysis including Dedekind cuts, limits and continuity.
\end{abstract}


\maketitle

\section{introduction}\label{sec:introduction}

Following the pioneering work of Bolzano and Weierstrass, ``$(\varepsilon,\delta)$-definitions'' are at the heart of textbook presentations of elementary analysis (e.g.\ \cite{rudin}). While with practice the motivated student quickly becomes proficient in this language, it is natural to ask if fundamental notions such as limits, continuity and integrals could perhaps be defined more conceptually. 

In the present paper we develop a rather general framework, which we refer to as {\it Darboux calculus}, whose specialization to the context of real analysis provides a unified and conceptual approach to all the main definitions encountered in, say, single variable calculus. Our starting point is the observation that the completeness of the ordered set of extended real numbers $\widehat \R=\{\pm \infty\}\cup \R$ is equivalent to the validity of the following

\begin{lem}\label{lem:intro}
Let $\O$ be a (partially) ordered set, let $\S\subseteq \O$ be any subset and let $\psi:\S\to \widehat \R$ be an order preserving function. Then the set of order preserving functions $f:\O\to \widehat \R$ whose restriction to $\S$ coincides with $\psi$ has a maximum and a minimum.
\end{lem}
In particular, such an order preserving function $\psi$ singles out a distinguished subset 
$\Dar(\psi)\subseteq \O$, the {\it Darboux set of $\psi$}, of elements on which the maximum and minimum extensions of $\psi$ coincide. Equivalently, $\Dar(\psi)$ can be thought of as the subset to which $\psi$ extends canonically. We denote this canonical extension by $\ex_\psi$. 

The prototypical example of this construction is provided by the Darboux integral. Let $\O$ denote the set of all bounded functions on an interval $[a,b]\subseteq \R$, let $\S$ be the subset of step functions and let $\psi$ be the function that to each step function assigns its integral defined naively in terms of signed areas of rectangles. In this case, as shown in Example \ref{ex:integrals} below, $\Dar(\psi)$ coincides with the set of Darboux integrable functions on $[a,b]$ and $\ex_\psi$ is the Darboux integral. 

This approach to the Darboux integral exemplifies the philosophy of this paper: naturally occurring pairs $(\mathcal X,\varphi)$ consisting of a class $\mathcal X$ of $\widehat \R$-valued functions and an order preserving function $\varphi:\mathcal X\to \widehat \R$ are of the form $(\Dar(\psi),\ex_\psi)$ for a suitable order preserving function $\psi$ defined on a subset $\S\subseteq \mathcal X$ of functions that ``obviously belong to $\mathcal X$''.

For instance, let $\O$ be the set of all sequences of real numbers, let $\S$ be the subset of sequences that are eventually constant and let $\psi$ be the function that to each sequence $\eta\in\S$ assigns the only value that $\eta$ attains infinitely many times. Then, as shown in as shown in Example \ref{ex:sequences} below, $\Dar(\psi)$ coincides with the set of convergent (possibly to $\pm \infty$) sequences and $\ex_\psi(f)=\lim_n f(n)$ for every $f\in \Dar(\psi)$.  The advantage here is that instead of having to come up with a clever $(\varepsilon,\delta)$-definition of limit of a sequence we only need to prescribe the obvious limit of an eventually constant sequence and the formalism of Darboux calculus automatically takes care of the general case. 

Similarly, let $\O$ be the set of all functions $f:\R\to \R$ and fix $x_0\in \R$. It is shown in  Example \ref{ex:continuous} that if $\S$ denotes the set of all functions that are constant on some open neighborhood of $x_0$ and $\psi$ is the function that to each $\eta\in \S$ assigns $\psi(\eta)=\eta(x_0)$, then $\Dar(\psi)$ is the set of functions that are continuous at $x_0$ and $\ex_\psi(f)=f(x_0)$ for all $f\in \Dar(\psi)$. Once again, given as only input the set of functions that are obviously 
continuous at $x_0$, our machinery returns the set of functions that are continuous at $x_0$ as output. We view this as an intuitive alternative to the standard $(\varepsilon,\delta)$-definition of continuity. 

The statement of Lemma \ref{lem:intro} holds more generally if $\widehat \R$ is replaced with any ordered set that is complete in the sense that every subset has a least upper bound and a greatest lower bound. Furthermore, the inclusion $\iota: \mathcal S \hookrightarrow\mathcal O$ can be replaced with an arbitrary embedding of ordered sets. In fact, the reader familiar with category theory will easily recognize the maximum and minimum extensions of $\psi$ in Lemma \ref{lem:intro} as, respectively, the left and right Kan extensions \cite{maclane} (assuming they exist) of $\psi$ along $\iota$. Similarly, the Darboux set of $\psi$ can be thought of as the equalizer of the left and right Kan extensions. Here we are implicitly using the standard interpretation of an ordered set $\O$ a category with objects are the elements $x\in \O$ and such that ${\rm Hom}(x,y)$ consists of a single element if $x\le y$ and is empty otherwise. From the vantage point of category theory, the present paper can be summarized as the observation that equalizers of left and right Kan extensions arise naturally in elementary analysis. While some of our propositions and theorems are particular instances of much more general results about left and right Kan extensions, we choose to give self-contained proofs in the case of ordered sets. In this way, we hope to provide evidence of the effectiveness of Darboux calculus as a stand-alone approach to the foundations of analysis that might be one day used to teach the subject at the undergraduate level.

An example of the flexibility of categorical thinking in this context comes from looking at the Yoneda embedding of an ordered set $\O$ into the set of order-preserving functions from $\O$ to the unique (up to a unique isomorphism) nontrivial ordered set with two elements. As it turns out, the Darboux set of the identity function of the image of the Yoneda embedding essentially coincides with the Dedekind-MacNeille completion of $\O$. While the idea of understanding Dedekind cuts in terms of presheaves is not new (see e.g.\ \cite{taylor}), our emphasis is on the fact that Darboux sets are not only effective in isolating interesting classes of $\R$-valued functions but can be used to construct $\R$ itself! In fact we show that with a little more effort, the field structure of $\R$ can also be recovered from that of $\Q$ in terms of Kan extensions. Our exposition appears to be somewhat more succinct, direct and self-contained than previous treatments of elementary analysis based on category theory e.g.\ \cite{univalent}, \cite{taylor-lambda}, \cite{edalat}. It would be interesting to carry out a detailed comparison between these approaches and the one presented here.

The paper is organized as follows. Section \ref{sec:preliminary} contains basic material on ordered sets and order preserving functions. In Section \ref{sec:darboux-sets} we introduce the main concepts used in this paper, including Darboux sets and Darboux extensions. Section \ref{sec:complete} is devoted to the notion of completeness defined here in terms of extensions of partially defined order preserving functions. As we show, our definition, which we refer to Darboux completeness, is in fact equivalent to the more familiar notion of Dedekind completeness. 
In Sections \ref{sec:CIC}-\ref{sec:completion} we discuss the Yoneda embedding and the Darboux completion of an arbitrary ordered set. In particular in Section \ref{sec:CIC} we use Darboux extensions to prove that completely integrally closed subgroups of automorphisms of a complete ordered set lift to automorphisms of the completion, a result that we use to construct the field operations on $\R$. Our strategy here can be thought of as a Darboux-theoretic version of the approach used in \cite{fuchs} to establish similar results directly at the level of Dedekind cuts. Once the real numbers are constructed, in Section \ref{sec:limits} we shift our attention to ordered sets of $\R$-valued functions. We prove that an $\R$-valued function $f$ has limit with respect to some filter basis $F$ (in the sense that each $\varepsilon$-neighborhood of the limit contains the image $f(\S)$ of some $\S\in F$) if and only if $f$ is in the Darboux set of the partial function defined by assigning to each function constant on some $\S\in F$ the only value that it attains on $\S$. This characterization of convergence with respect to a filter basis yields at once Darboux-theoretic formulations of several $(\epsilon,\delta)$-definitions such as limits of sequences, limits of functions of one real variable and continuity. After discussing Darboux integrability (after which the general notion of Darboux set is modeled), we use Darboux calculus to prove a theorem which simultaneously generalizes the usual linearity theorems for limits, continuous functions and integrals. In fact, all the major theorems of elementary real analysis (e.g.\ the Intermediate Value Theorem, the Extreme Value Theorem and the Fundamental Theorem of Calculus) can be proved conceptually using the language of Darboux calculus.  We hope to come back to this point elsewhere and ultimately provide an exhaustive and fully self-contained treatment of elementary real analysis in the language of this paper.

\section{Preliminaries on Ordered Sets}\label{sec:preliminary}

\begin{dfn}
A {\it (partially) ordered set} is a set $\O$ together with a reflexive, antisymmetric, and transitive relation which we denote by $\le$.
\end{dfn}

\begin{example}\label{ex:2.2}
If $\O$ is an ordered set, every subset $\S\subseteq \O$ inherits an induced order. For every $x,y\in \O$ such that $x\le y$, the {\it interval with endpoints $x$ and $y$} is the (ordered) subset $[x,y]$ of all $z\in \O$ such that $x\le z\le y$. 
\end{example}

\begin{example}
A {\it discrete} set is an ordered set with the trivial order with respect to which $x\le y$ if and only if $x=y$. If $\mathcal O$ is an ordered set, we denote by $|\O|$ its underlying discrete set. 
\end{example}

\begin{rem}
If $\O$ is an ordered set, we denote by $\O^\op$ the {\it opposite} ordered set such that $|\O^{\rm op}|=|\O|$ and $x\le y$ in $\O^\op$ if and only if $y\le x$ in $\mathcal O$. 
\end{rem}

\begin{example}
Given two ordered sets $\O_1,\O_2$, we denote by $\O_1\times \O_2$ the ordered set such that $|\O_1\times \O_2|=|\O_1|\times |\O_2|$ with order such that $(x_1,x_2)\le (y_1,y_2)$ if and only if $x_1\le y_1$ and $x_2\le y_2$.
\end{example}

\begin{dfn}
Let $\O$ and $\P$ be ordered sets. The set of {\it order preserving functions from $\O$ to $\P$} is
\[
\OP(\O,\P)=\{f:|\O|\to |\P|\,|\,f(x)\le f(y) \text{ if } x\le y\}\,. 
\]
We view $\OP(\O,\P)$ as an ordered set such that $f\le g$ if and only if $f(x)\le g(x)$ for all $x\in \O$. We use the shorthand notation $\OP(\O)=\OP(\O,\O)$. We also say that $f\in \OP(\O,\P)$ is an {\it embedding} if for any $x,y\in \mathcal O$, $f(x)\le f(y)$ implies $x\le y$. An {\it isomorphism} is a surjective embedding. Given an ordered set $\O$, we denote by $\Aut(\O)$ the group of all isomorphisms in $\OP(\O)$. 
\end{dfn}

\begin{dfn}
If $\O$ is an ordered set, we define its {\it augmentation} to be the ordered set $\widehat \O$ such that 
\begin{enumerate}[1)]
\item $|\widehat \O|=|\O|\cup\{-\infty,+\infty\}$;
\item the canonical inclusion of $|\O|$ into $|\widehat \O|$ defines an embedding of $\O$ into $\widehat \O$;
\item $\widehat \O=[-\infty,+\infty]$.
\end{enumerate}
\end{dfn}

\begin{dfn}
Let $\O$ and $\P$ be ordered sets. A {\it partial function $\psi:\O\hto \P$ from $\O$ to $\P$} is an order preserving function $\psi:\Dom(\psi)\to \P$ defined on an ordered subset $\Dom(\psi)\subseteq \O$ called the {\it domain} of $\psi$. The ordered set $\im(\psi)=\psi(\Dom(\psi))$ is called the {\it image} of $\psi$. An {\it extension of $\psi$ to $\O$} is an order preserving function $f:\O\to \P$ whose restriction $f_{|_\Dom(\psi)}$ to $\Dom(\psi)$ coincides with $\psi$. 
\end{dfn}

\begin{example}\label{ex:2.9}
Let $\O$ be an ordered set and let $\1$ be the unique (up to a unique isomorphism) ordered set with one element. Then $\O$ is canonically identified with $\OP(\1,\O)$.
\end{example}

\begin{dfn} Let $\O$ and $\P$ be ordered sets. A set $\Psi$ of partial functions from $\O$ to $\P$ is {\it compatible} if for any $\psi',\psi''\in \Psi$, the restrictions of $\psi'$ and $\psi''$ to $\Dom(\psi')\cap \Dom(\psi'')$ coincide. If $\Psi$ is compatible, we define its {\it common extension} to be the partial function $\psi:\O\hto\P$ such that 
\[
\Dom(\psi)=\bigcup_{\psi'\in\Psi} \Dom(\psi')
\]
and $\psi(x)=\psi'(x)$ for every $x\in \Dom(\psi')$ and for every $\psi'\in \Psi$.
\end{dfn}

\begin{rem}\label{rem:partition-1}Let $\O$ and $\P$ be ordered sets.
If $\psi$ is the common extension of a compatible set $\Psi$ of partial functions from $\O$ to $\P$, then $f:\O\to P$ is an extension of $\psi$ to $\O$ if and only if it is an extension of $\psi'$ to $\O$ for each $\psi'\in \Psi$.
\end{rem}

\section{Darboux Sets and Darboux Extensions}\label{sec:darboux-sets}

\begin{dfn} Let $\O$ and $\P$ be ordered sets.
A partial function $\psi:\O\hto \P$ is {\it extremizable} if there exist order preserving functions $\lex_\psi,\uex_\psi:\O\to\P$ such that
$\lex_\psi\le f\le \uex_\psi$ for all extensions $f:\O\to \P$ of $\psi$ to $\O$.
If this is the case, we call $\lex_\psi$ (resp.\ $\uex_\psi$) the {\it lower extension} (resp.\ the {\it upper extension}) of $\psi$.
\end{dfn}

\begin{rem}
Let $\psi:\O\to \P$ be an extremizable partial function. If $x\in \Dom(\psi)$, then $\lex_\psi(x)=\uex_\psi(x)$. Therefore, if $f:\O\to \P$ is an order preserving function such that $\lex_\psi\le f\le \uex_\psi$, then $f$ is automatically an extension of $\psi$. 
\end{rem}

\begin{example}\label{ex:deltafunction}
Let $\O$ and $\P$ be ordered sets and let $\psi:\O\hto\widehat \P$ be a partial function such that $\Dom(\psi)=\{x\}$. Then $\psi$ is extremizable. Moreover, $\lex_\psi(y)=\psi(x)$ if $x\le y$ and $\lex_\psi(y)=-\infty$ otherwise. Similarly,  $\uex_\psi(y)=\psi(x)$ if $y\le x$ and $\uex_\psi(y)=+\infty$ otherwise. 
\end{example}

\begin{dfn} Let $\O$ and $\P$ be ordered sets. For each extremizable partial function $\psi:\O\hto\P$, we define the {\it Darboux set of $\psi$} to be
\[
\Dar(\psi)=\{x\in \O\,|\, \lex_\psi(x)=\uex_\psi(x)\}\,.
\]
Moreover, we denote by $\ex_\psi:\O\hto\P$ the {\it Darboux extension of $\psi$}, i.e.\ the restriction of $\uex_\psi$ (or equivalently of $\lex_\psi$) to $\Dar(\psi)$. 
\end{dfn}

\begin{dfn} Let $\O$, $\P$ be ordered sets and let $\psi$ be a partial function from $\O$ to $\P$. We say that $x\in \O$ is {\it $\psi$-bounded} if $y\le x\le z$ for some $y,z\in \Dom(\psi)$. We denote the set of $\psi$-bounded elements of $\O$ by $\B(\psi)$. We say that $\psi$ is {\it encompassing} if every element of $\O$ is $\psi$-bounded.
Moreover, for each extremizable $\psi:\O\hto\P$ we define the {\it bounded Darboux set of $\psi$} to be
the subset $\BDar(\psi)$ of all $\psi$-bounded elements of $\Dar(\psi)$. 
\end{dfn}

\begin{rem}\label{rem:partition-2} Let $\O$, $\P$ be ordered sets and let $\psi$ be the common extension of a compatible set $\Psi$ of partial functions from $\O$ to $\P$. If any $\psi'\in \Psi$ is encompassing, then $\Dom(\psi')\subseteq \Dom(\psi)$ implies that $\psi$ is also encompassing.
\end{rem}

\begin{rem}\label{rem:partition-3} Let $\O$ and $\P$ be ordered sets and let $\Psi$ be a compatible set of extremizable partial functions from $\O$ to $\P$. If the common extension $\psi$ of $\Psi$ is also extremizable, then Remark \ref{rem:partition-1} implies that $f\in [\lex_\psi,\uex_\psi]$ if and only if $f\in [\lex_{\psi'},\uex_{\psi'}]$ for each $\psi'\in \Psi$. In particular,
\[
\{\uex_\psi\}=\bigcap_{\psi'\in \Psi} [\uex_\psi,\uex_{\psi'}]\qquad\textrm{and}\qquad
\{\lex_\psi\}=\bigcap_{\psi'\in \Psi} [\lex_{\psi'},\lex_\psi]\,.
\]
\end{rem}

\begin{rem} Let $\O$ and $\P$ be ordered sets and let 
$\psi$ be an extremizable partial function from $\O$ to $\P$. If $f:\O\to \P$ is an extension of $\ex_\psi$ to $\O$, then its restriction to $\Dom(\psi)$ coincides with $\psi$ 
and thus $\lex_\psi\le f\le \uex_\psi$. Since by construction $\lex_\psi$ and $\uex_\psi$ restrict to $\ex_\psi$ on $\Dar(\psi)$, it follows that the set of extensions of $\ex_\psi$ to $\O$ coincides with the set of extension of $\psi$ to $\O$. In particular, $\Dar(\ex_\psi)=\Dar(\psi)$ and $\ex_{\ex_\psi}=\ex_\psi$.
\end{rem}

\begin{dfn}
Let $\O_1$, $\O_2$ and $\O_3$ be ordered sets. The partial functions $\psi_1:\O_1\hto\O_2$ and $\psi_2:\O_2\hto\O_3$ are {\it composable} if $\Dom(\psi_2)\cap \im(\psi_1)$ is non-empty. If this is the case, their {\it composition} is the partial function $\psi_2\circ\psi_1:\O_1\hto \O_3$ such that $(\psi_2\circ\psi_1)(x)=\psi_2(\psi_1(x))$ for each $x$ in
\[
\Dom(\psi_2\circ \psi_1)=\{x\in \Dom(\psi_1)\,|\,\psi_1(x)\in \Dom(\psi_2)\}\,.
\]

\end{dfn}

\begin{prop}\label{prop:inequality}
Let $\O_1$, $\O_2$ and $\O_3$ be ordered sets. Let $\psi_1:\O_1\hto\O_2$, $\psi_2:\O_2\hto\O_3$ be partial functions that
\begin{enumerate}[i)]
\item $\Dom(\psi_2)\subseteq \im(\psi_1)$;
\item $\psi_1$, $\psi_2$ and $\psi_2\circ\psi_1$ are extremizable.
\end{enumerate}
Then
\begin{enumerate}[1)]
\item $\lex_{\psi_2\circ\psi_1} \le \lex_{\psi_2}\circ \lex_{\psi_1}\le\uex_{\psi_2}\circ \uex_{\psi_1}\le \uex_{\psi_2\circ \psi_1}$;
\item $\ex_{\psi_1}(\Dar(\psi_2\circ\psi_1)\cap \Dar(\psi_1))\subseteq \Dar(\psi_2)$;
\item $(\ex_{\psi_2}\circ \ex_{\psi_1})(x)=\ex_{\psi_2\circ \psi_1}(x)$ for all $x\in \Dar(\psi_1)\cap \Dar(\psi_2\circ \psi_1)$.
\end{enumerate}
\end{prop}

\proof 1) is a consequence of the fact that $\uex_{\psi_2}\circ \uex_{\psi_1}$ and $\lex_{\psi_2}\circ \lex_{\psi_1}$ are extensions of $\psi_2\circ \psi_1$ to $\O_1$. If $x\in \Dar(\psi_2\circ \psi_1)\cap \Dar(\psi_1)$, then 1) implies
\[
\ex_{\psi_2\circ \psi_1}(x)=\lex_{\psi_2}( \ex_{\psi_1}(x))=\uex_{\psi_2}(\ex_{\psi_1}(x))
\]
which proves 2) and 3). 

\begin{rem}\label{rem:inequality} Since $\Dom(\psi_2)=\O_2$ implies $\lex_{\psi_2}=\psi_2=\uex_{\psi_2}$, then $\Dar(\psi_2\circ \psi_1)\subseteq \Dar(\psi_1)$ whenever the partial function $\psi_2$ in the statement of Proposition \ref{prop:inequality} is an embedding and thus $\ex_{\psi_1}(\Dar(\psi_2\circ \psi_1))\subseteq \Dar(\psi_2)$.
\end{rem}

\begin{lem}\label{lem:ev}
Let $\O$, $\P$, $\P'$ be ordered sets and let $\psi:\O\hto\OP(\P,\P')$ be an extremizable partial function. For each $p\in \P$, let $\ev_p:\OP(\P,\P')\to \P'$ be the order preserving function that to each $f:\P\to \P'$ assigns its evaluation $\ev_p(f)=f(p)$ at $p$. If $\ev_p\circ \psi:\O\hto \P'$ is extremizable for every $p\in \P$, then
\[
\ev_p\circ \uex_\psi =\uex_{\ev_p\circ \psi} \quad \textrm{ and } \quad \ev_p\circ \lex_\psi = \lex_{\ev_p\circ \psi}\,.
\]
\end{lem}

\proof Using Proposition \ref{prop:inequality}, $\ev_p\circ \uex_\psi = \uex_{\ev_p} \circ \uex_\psi \le \uex_{\ev_p\circ \psi}$. Consider the order preserving function $g:\O\to\OP(\P,\P')$ such that $(g(x))(p)=\uex_{\ev_p\circ \psi}$ for every $x\in \O$ and for every $p\in \P$. Then $(g(\eta))(p)=\uex_{\ev_p\circ \psi}(\eta) = (\psi(\eta))(p)$ for every $\eta \in \Dom(\psi)$. Therefore, $g\le \uex_\psi$ and thus $\uex_{\ev_p\circ \psi} = \ev_p\circ g \le \ev_p\circ \uex_\psi$. Hence, $\ev_p\circ \uex_\psi = \uex_{\ev_p\circ \psi}$. The second equality is proved in a similar way.

\begin{lem}\label{lem:product}
Let $\O$, $\P_1$, $\P_2$ be ordered sets, let $\psi:\O\hto \P_1\times \P_2$ be a partial function and and for $i=1,2$ let $\pi_i:\P_1\times \P_2\to \P_i$ be the (order preserving) projection onto the respective factor. Then $\psi$ is extremizable if and only if $\pi_i\circ \psi:\O\hto \P_i$ is extremizable for each $i=1,2$. If this is the case, then $\pi_i\circ \uex_\psi=\uex_{\pi_i\circ \psi}$ and $\pi_i\circ \lex_\psi=\lex_{\pi_i\circ \psi}$ for each $i=1,2$. 
\end{lem}

\proof Assume that $\psi$ is extremizable. Then $\pi_i\circ \lex_\psi$ and $\pi_i\circ \uex_\psi$ are extensions of the partial function $\pi_i\circ \psi:\O\hto \P_i$ (with domain $\Dom(\psi)$) for each $i=1,2$. Furthermore, if $f_1:\O\to \P_1$ and $f_2:\O\to \P_2$ are, respectively, extensions of $\pi_1\circ \psi$ and $\pi_2\circ \psi$, then $(f_1,f_2):\O\to \P_1\times \P_2$ is an extension of $\psi$. By assumption, this implies $\lex_\psi\le (f_1,f_2)\le \uex_\psi$ and thus
\[
\pi_i\circ\lex_\psi \le f_i\le \pi_i\circ \uex_\psi 
\]
for each $i=1,2$. Hence $\pi\circ \psi_i$ is extremizable, $\pi_i\circ \uex_\psi=\uex_{\pi_i\circ \psi}$ and $\pi_i\circ \lex_\psi=\lex_{\pi_i\circ \psi}$ for each $i=1,2$. Conversely, assume that $\pi_1\circ \psi$ and $\pi_2\circ\psi$ are extremizable. Then $(\lex_{\pi_1\circ\psi},\lex_{\pi_2\circ \psi})$
and $(\uex_{\pi_1\circ\psi},\uex_{\pi_2\circ \psi})$ are both extensions of $\psi=(\pi_1\circ\psi,\pi_2\circ \psi)$. Moreover, if $f:\O\to \P_1\times \P_2$ is any extension of $\psi$, then $\pi_i\circ f$ is an extension of $\pi_i\circ \psi$ for each $i=1,2$. Since $f=(\pi_1\circ f,\pi_2\circ f)$, this implies 
\[
(\lex_{\pi_1\circ\psi},\lex_{\pi_2\circ \psi})\le f \le (\uex_{\pi_1\circ\psi},\uex_{\pi_2\circ \psi})
\]
and thus $\psi$ is extremizable with lower extension equal to $(\lex_{\pi_1\circ\psi},\lex_{\pi_2\circ \psi})$ and upper extension equal to $(\uex_{\pi_1\circ\psi},\uex_{\pi_2\circ \psi})$.

\begin{rem}\label{rem:3.12}
Let $\O$ be a non-empty ordered set and let $\emptyset:\O\hto \O$ be {\it the empty partial function of $\O$} i.e.\ the unique partial function from $\O$ to itself whose domain is the empty set. 
Since the set of extensions of $\emptyset$ to $\O$ coincides with $\OP(\O)$, if $\emptyset$ is extremizable, then in particular $\lex_\emptyset\le x\le \uex_\emptyset$ for every constant function $x:\O\to \O$. In other words, the lower and upper Darboux extensions of the empty partial function are constant and (with a slight abuse of notation) $\O=[\lex_\emptyset(\O),\uex_\emptyset(\O)]$.
\end{rem}

\section{Darboux complete ordered sets}\label{sec:complete}
\begin{dfn}
An ordered set $\P$ is {\it Darboux complete} if every partial function from $\widehat \P$ to itself is extremizable.
\end{dfn}

\begin{example}\label{ex:4.2}
Since by Example \ref{ex:deltafunction} each partial function $f:\widehat\emptyset\hto\widehat\emptyset$ is extremizable, the empty ordered set $\emptyset$ is Darboux complete.
\end{example}

\begin{lem}\label{lem:complete}
Let $\S$ be a non-empty subset of a Darboux complete ordered set $\P$. If $\id_\S:\widehat\P\hto\widehat\P$ denotes the identity function on $\S$ and $J=[\uex_{\id_\S}(-\infty),\lex_{\id_\S}(+\infty)]$, then
\begin{enumerate}[1)]
\item $\S\subseteq J$;
\item $J$ is the intersection of all intervals of $\widehat \P$ that contain $\S$.
\end{enumerate}
\end{lem}

\proof Since $\P$ is Darboux complete, then $\lex_{\id_\S}$ and $\uex_{\id_\S}$ exist. For every $s\in \S$
\[
\uex_{\id_\S}(-\infty)\le \uex_{\id_\S}(s)=\lex_{\id_\S}(s)\le \lex_{\id_\S}(+\infty)\,,
\]
which implies 1). If $x,y\in\widehat\P$ are such that $\S\subseteq [x,y]$, let $\psi:\widehat \P\hto\widehat\P$ be the partial function with domain $\widehat \S$, whose restriction to $\S$ is the identity and such that $\psi(-\infty)=x$ and $\psi(+\infty)=y$. Then
\[
x=\uex_\psi(-\infty)\le \uex_{\id_\S}(-\infty) \le \lex_{\id_\S}(+\infty)\le \lex_\psi(-\infty)=y\,,
\]
which concludes the proof.

\begin{prop}\label{prop:4.4}
Let $\P$ be an ordered set. The following are equivalent
\begin{enumerate}[1)]
\item $\P$ is Darboux complete;
\item Every partial function with codomain $\widehat \P$ is extremizable; 
\item for every ordered set $\O$, every partial function with codomain $\OP(\O,\widehat \P)$ is extremizable.
\end{enumerate}
\end{prop}

\proof Assume that $\P$ is Darboux complete. Let $\psi$ be a partial function from an ordered set $\O$ to $\widehat\P$ and let $x\in \O$. Consider the subsets
\begin{equation}\label{eq:S_x}
\S_x=\{\psi(y)\,|\,y\le x \textrm{ and } y\in\Dom(\psi)\}\subseteq \widehat \P
\end{equation}
and
\begin{equation}\label{eq:S^x}
\S^x=\{\psi(y)\,|\,x\le y \textrm{ and } y\in\Dom(\psi)\}\subseteq \widehat \P
\end{equation}
together with their identity functions $\id_{\S_x},\id_{\S^x}:\widehat \P\hto\widehat\P$. Define $l,u:\O\to\widehat P$ such that
\[
l(x)=\lex_{\id_{\S_x}}(+\infty)\qquad \textrm{and}\qquad u(x)=\uex_{\id_{\S^x}}(-\infty)
\]
for all $x\in \O$. To see that $l$ and $u$ are indeed order preserving, assume that $x,y\in \O$ are such that $x\le y$. Since $\S_x\subseteq\S_y$, then $\lex_{\id_{\S_y}}$ is an extension of $\id_{\S_x}$ and thus $l$ is order preserving. Similarly, $u$ is order preserving because $\S^y\subseteq \S^x$ implies that the restriction of $\lex_{\id_{\S_{x}}}$ to $\S^y$ coincides with $\id_{\S^y}$. Moreover $l$ is an extension of $\psi$ to $\O$ since for every $x\in \Dom(\psi)$, $\S_x\subseteq [-\infty,\psi(x)]$ and Lemma \ref{lem:complete} implies that
\[
\psi(x)=\lex_{\id_{\S_x}}(\psi(x))\le l(x)\le \psi(x)\,.
\]
On the other hand, $\psi(y)=f(y)\le f(x)$ for any extension $f$ of $\psi$ to $\O$ and for any $\psi(y)\in \S_x$. Therefore, $\S_x\subseteq [-\infty, f(x)]$ and thus (using again Lemma \ref{lem:complete}), $l(x)\le f(x)$ . Together with a similar argument involving $u$, this proves 2). Assume that 2) holds and let $\O$, $\O'$ be arbitrary ordered sets. Consider the canonical embedding $\alpha$ that to each partial function $\psi:\O'\hto\OP(\O,\widehat \P)$ assigns the partial function $\alpha(\psi):\O'\times \O\hto\widehat\P$ such that $(\alpha(\psi))(x',x)=(\psi(x'))(x)$ for all $(x',x)\in \Dom(\alpha(\psi))=\Dom(\psi)\times \O$. The Darboux completeness of $\P$ ensures that $\alpha(\psi)$ is extremizable and thus $\lex_{\alpha(\psi)}\le \alpha(f)\le \uex_{\alpha(\psi)}$ for each extension $f$ of $\psi$ to $\O'$. Since the restriction of $\alpha$ to the subset of order preserving functions $\O'\to\OP(\O,\widehat \P)$ is an isomorphism, then $\lex_\psi=\alpha^{-1}(\lex_{\alpha(\psi)})$ and $\uex_\psi=\alpha^{-1}(\uex_{\alpha(\psi)})$, which proves 3). Example \ref{ex:2.9} shows that 1) is a particular case of 3), which concludes the proof. 

\begin{rem}\label{rem:4.9}
Let $\P$ be an ordered set. Assume $\P$ is a Darboux complete ordered set, and $\S\subseteq \P$ is nonempty and bounded i.e.\ $\S\subseteq [x,y]$ for some $x,y\in \P$. Then Lemma \ref{lem:complete} implies that $\lex_{\id_\S}(+\infty)$ and $\uex_{\id_\S}(-\infty)$ are respectively the least upper bound $\sup(\S)$ and the greatest lower bound $\inf(\S)$ of $\S$. Therefore, $\P$ is Dedekind complete. Conversely, suppose that the least upper bound and the greatest lower bound of every nonempty bounded subset of $\P$ exist. Given any partial function $\psi:\widehat \P\to \widehat\P$, let $\S_x$ and $\S^x$ be defined as in \eqref{eq:S_x} and \eqref{eq:S^x} respectively. Then the same argument as in the proof of Proposition \ref{prop:4.4} shows that $\psi$ is extremizable with $\lex_\psi(x)=\sup(\S_x)$ and $\uex_\psi(x)=\inf(\S^x)$ for all $x\in \O$. Hence, $\P$ is Darboux complete if and only if $\P$ is Dedekind complete. While these two notions of completeness are equivalent, the point of view of this paper is that Darboux completeness allows for a more direct and conceptual route to the foundations of elementary analysis.  
\end{rem}

\begin{cor}\label{cor:encompassing}
Let $\O$ be an ordered set, let $\P$ be a Darboux complete ordered set and let $N$ be a positive integer. Every encompassing partial function from $\O$ to $\P^N$ is extremizable. 
\end{cor}

\proof By Lemma \ref{lem:product}, it suffices to prove the $N=1$ case. Let $\varphi:\O\hto \P$ be encompassing and let $\iota:\P\to \widehat \P$. By assumption, for each $z\in \O$ there exist $x,y\in \Dom(\varphi)$ such that $x\le z\le y$ and thus 
\[
\varphi(x)=\iota(\varphi(x))\le \lex_{\iota\circ\varphi}(z)\le \uex_{\iota\circ\varphi}(z)\le \iota(\psi(y))=\varphi(y)\,.
\]
Therefore, $\lex_{\iota\circ\varphi}$ and $\uex_{\iota\circ\varphi}$ have their image contained in $\P$ and thus  are extensions of $\varphi$ to $\O$. Moreover, $\lex_{\iota\circ\varphi}(x)\le f(x)\le \uex_{\iota\circ\varphi}(x)$ for every extension $f$ of $\varphi$ to $\O$ and for every $x\in \P$. Hence $\varphi$ is extremizable and $\lex_\varphi(x)=\lex_{\iota\circ\varphi}(x)$, $\uex_\varphi(x)=\uex_{\iota\varphi}(x)$ for all $x\in \O$.

\begin{example}\label{ex:4.6}
We define the {\it free cocompletion} of an ordered set $\O$ to be the ordered set $\O^\vee = \OP(\O^\op,\widehat\emptyset)$. Let $\P=\O^\vee\setminus\{\pm\infty\}$, where $\pm\infty$ denotes the constant function such that $\im(\pm\infty)=\pm\infty$. Combining Example \ref{ex:4.2} with Proposition \ref{prop:4.4} shows that every partial function with codomain $\O^\vee=\widehat{\P}$ is extremizable. Using Proposition \ref{prop:4.4} again, we conclude that $\P$ is Darboux complete. 
\end{example}

\begin{example}\label{ex:constant}
Let $\O$ be an ordered set, let $\P$ a Darboux complete ordered set and let $\S$ be a nonempty subset of $\O$. Furthermore, let $\psi_\S:\OP(\O,\P)\hto\widehat \P$ be the partial function with domain the subset of functions that are constant on $\S$ and such that $\psi_\S(f)=f(x)$ for every $f\in \Dom(\psi_\S)$ and every $x\in \S$. For each $x\in \S$, $\ev_x$ coincides with $\psi_\S$ on $\Dom(\psi_\S)$ and thus $\ev_x\in[\lex_{\psi_\S},\uex_{\psi_\S}]$. In particular, if $f:\O\to \P$ is in the Darboux set of $\psi_\S$, then $\ev_x\circ f=\ev_y\circ f$ for every $x,y\in \S$ i.e.\ $f$ is constant on $\S$. Hence, $\Dom(\psi_\S)=\Dar(\psi_\S)$. 
\end{example}

\begin{rem}\label{rem:constant} Using the notation of Example \ref{ex:constant}, assume furthermore that $\O$ is discrete. For every order preserving function $f:\O\to\P$ and for every $y\in \P$, let $f_y\in \Dom(\psi_\S)$ be the function whose restriction to $\O\setminus \S$ coincides with $f$ and such that $f_y(x)=y$ for all $x\in \S$. In particular, if there exists $y,z\in \P$ such that $f(x)\in [y,z]$ for all $x\in \S$, then $f\in [f_y,f_z]$ and thus $[\lex_{\psi_\S},\uex_{\psi_\S}]\subseteq [y,z]$. Moreover, Corollary \ref{cor:encompassing} implies that the restriction $\varphi_\S:\B(\psi_\S)\hto\P$ of $\psi_\S$ to $\B(\psi)$ is extremizable. 
\end{rem}

\begin{rem}\label{rem:4.9}
Let $\P$ be an ordered set. Assume $\P$ is a Darboux complete ordered set, and $\S\subseteq \P$ is nonempty and bounded i.e.\ $\S\subseteq [x,y]$ for some $x,y\in \P$. Then Lemma \ref{lem:complete} implies that $\lex_{\id_\S}(+\infty)$ and $\uex_{\id_\S}(-\infty)$ are respectively the least upper bound $\sup(\S)$ and the greatest lower bound $\inf(\S)$ of $\S$. Therefore, $\P$ is Dedekind complete. Conversely, suppose that the least upper bound and the greatest lower bound of every nonempty bounded subset of $\P$ exist. Given any partial function $\psi:\widehat \P\to \widehat\P$, let $\S_x$ and $\S^x$ be defined as in \eqref{eq:S_x} and \eqref{eq:S^x} respectively. Then the same argument as in the proof of Proposition \ref{prop:4.4} shows that $\psi$ is extremizable with $\lex_\psi(x)=\sup(\S_x)$ and $\uex_\psi(x)=\inf(\S^x)$ for all $x\in \O$. Hence, $\P$ is Darboux complete if and only if $\P$ is Dedekind complete. While these two notions of completeness are equivalent, the point of view of this paper is that Darboux completeness allows for a more direct and conceptual route to the foundations of elementary analysis.  
\end{rem}

\section{Completely Integrally Closed Subgroups}\label{sec:CIC}

\begin{prop}\label{prop:product}
Let $\O$, $\O'$ be ordered sets, let $\P$ be a Darboux complete ordered set and consider the composition of ordered function $\mu:\OP(\widehat \P)\times\OP(\widehat\P)\to\OP(\widehat \P)$ defined by setting $\mu(\varphi,\varphi')=\varphi\circ\varphi'$ for all $\varphi,\varphi'\in \OP(\widehat\P)$. If $\psi:\O\hto\OP(\widehat \P)$ and $\psi':\O'\hto\OP(\widehat \P)$ are partial functions with images in $\Aut(\widehat\P)$, then
\[
\mu\circ(\uex_\psi \times \uex_{\psi'}) = \uex_{\mu\circ(\psi\times \psi')}\quad \textrm{ and }\quad \mu\circ (\lex_\psi \times \lex_{\psi'}) = \lex_{\mu\circ(\psi\times \psi')}\,.
\]

\end{prop}

\proof Since $\mu\circ (\uex_\psi \times \uex_{\psi'})$ is an extension of $\mu\circ (\psi\times \psi')$ to $\O\times \O'$, then $\mu\circ (\uex_\psi\times \uex_{\psi'})\le \uex_{\mu\circ (\psi\times \psi')}$.
On the other hand, if $\eta\in \Dom(\psi)$ is fixed, then 
\[
(\psi(\eta))^{-1}\circ\uex_{\mu \circ (\psi\times \psi')}(\eta,\eta')=\psi'(\eta')
\]
for every $\eta'\in \Dom(\psi')$. Using the assumption that $\psi'$ is extremizable, it follows that 
\[
(\psi(\eta))^{-1}\circ\uex_{\mu\circ(\psi\times \psi')}(\eta,x')\le \uex_{\psi'}(x')
\]
and thus
\[
\uex_{\mu\circ(\psi\times \psi')}(\eta,x')\le \psi(\eta)\circ\uex_{\psi'}(x')=(\mu \circ (\uex_\psi\times \uex_{\psi'}))(\eta,x') \le \uex_{\mu\circ(\psi\times \psi')}(\eta,x')
\]
for all $(\eta,x')\in \Dom(\psi) \times \O'$. Setting $q=\left(\uex_{\psi'}(x')\right)(p)$ yields
\[
\ev_p\circ \uex_{\mu\circ (\psi\times \psi')}(\eta,x)=(\psi(\eta))(q)=\uex_{\ev_q\circ\psi} (\eta)
\]
for every $p\in \widehat\P$ and for every $\eta\in \Dom(\psi)$. Lemma \ref{lem:ev} then implies
\[
\ev_p\circ \uex_{\mu\circ(\psi\times \psi')}(x,x') \le \uex_{\ev_q \circ \psi}(x) = \ev_q \circ \uex_\psi(x) = \ev_p\circ \mu\circ (\uex_\psi \times \uex_{\psi'})(x,x')
\]
for all $p\in \widehat \P$ and for all $(x,x')\in \O\times \O'$. This proves the first half of the Proposition, the second equality is proved in a similar way.

\begin{rem}
Given any ordered set $\O$, the ordered set $\OP(\O)$ of order preserving functions $f:\O\to \O$ is a monoid with respect to composition.
\end{rem}

\begin{dfn}
Let $\O$ be an ordered set. A subgroup (that is a sub-monoid closed under inverses) $\A$ of $\OP(\O)$ is {\it completely integrally closed} if for every $a,a'\in \A$, $a^n\le a'$ for all $n\in \N$ implies $a\le \id_\O$. 
\end{dfn}

\begin{rem} Completely integrally closed subgroups are a particular instance of the more general notion of (abstract) completely integrally closed ordered groups which plays a key role in the classical study \cite{fuchs} of embeddings in Dedekind complete ordered groups. The reminder of this section can be thought of as an alternate construction of these embeddings formulated in the equivalent language of Darboux complete ordered sets. Our main application is the self-contained construction of the field structure on the ordered set of real numbers described in Section \ref{sec:completion}.
\end{rem}

\begin{prop}\label{prop:CIC}
Let $\P$ be a Darboux complete ordered set. If $\A$ is a completely integrally closed subgroup of $\OP(\widehat\P)$, then $\BDar(\id_\A)$ is a subgroup of $\OP(\widehat\P)$. 
\end{prop}

\proof Since $\lex_{\id_\A}\circ \mu$  and $\uex_{\id_\A}\circ \mu$ to $\A\times \A$ are extensions of $\mu\circ (\id_\A\times \id_\A)$ to $(\OP(\widehat \P))^2$, we obtain
 \begin{equation}\label{eq:CIC1}
\lex_{\mu\circ (\id_\A \times \id_\A)} \le \lex_{\id_\A}\circ \mu\le \uex_{\id_\A}\circ \mu\le \uex_{\mu\circ (\id_\A \times \id_\A)}\,.
\end{equation}
By Proposition \ref{prop:product}, we conclude that these inequalities restrict to equalities on $(\Dar(\id_\A))^2$. Hence $\Dar(\id_\A)$ is closed under composition. We conclude that $\Dar(\id_\A)$, which contains the submonoid $\A$ of $\OP(\widehat \P)$, is itself a submonoid of $\OP(\widehat\P)$. Given $\varphi_1,\varphi_2\in \BDar(\id_\A)$, by definition there exist $a_i,a_i'\in \A$ such that $a_i\le \varphi_i\le a_i'$ for $i=1,2$. Therefore $a_1\circ a_2\le \varphi_1\circ \varphi_2\le a_1'\circ a_2'$ and thus $\BDar(\psi)$ is also a submonoid. In order to construct inverses, consider the partial function $\psi:(\OP(\widehat\P))^\op\hto \OP(\widehat\P)$ with domain $\A$ and such that $\psi(a)=a^{-1}$ for every $a\in \A$. By Proposition \ref{prop:product},
\begin{equation*}
\lex_\psi(\varphi) \circ \varphi=\lex_{\mu \circ (\psi\times \id_\A)}(\varphi,\varphi)\le \uex_{\mu \circ (\psi\times \id_\A)}(\varphi,\varphi) 
\end{equation*}
for all $\varphi\in \Dar(\id_\A)$. Since $\im(\psi)=\A$, then
\[
\id_\A\circ \mu(\psi\times \id_\A) = \mu(\psi\times \id_\A)
\]
and thus, using Proposition \ref{prop:inequality},
\begin{equation}\label{eq:CIC2}
\lex_\psi(\varphi) \circ \varphi\le \lex_{\id_\A}\left(\lex_{\mu \circ (\psi\times \id_\A)}(\varphi,\varphi)\right)\le \lex_{\id_\A}\left(\uex_{\mu \circ (\psi\times \id_\A)}(\varphi,\varphi)\right)\,. 
\end{equation}
Given $\varphi\in \BDar(\id_\A)$, $a\in \A$ such that $a\le \varphi$, and $a'\in \A$ such that $a'\le \uex_{\mu\circ(\psi\times \id_\A)}(\varphi,\varphi)$, then 
\[
a\circ a'\le a \circ \uex_{\mu\circ(\psi\times \id_\A)}(\varphi,\varphi) \le a \circ \uex_{\mu\circ(\psi\times \id_\A)}(a,\varphi) = \varphi\,.
\]
Iterating the same argument with $a$ replaced by $a\circ (a')^{n-1}$, yields $a\circ (a')^n\le \varphi$ for all $n\in \N$. Since $\A$ is completely integrally closed, this implies $a'\le \id_{\widehat\P}$. Together with a similar argument involving $\lex_{\mu\circ(\psi\times \id_\A)}$, we conclude that 
\[
\lex_{\id_\A}\left(\uex_{\mu \circ (\psi\times \id_\A)}(\varphi,\varphi)\right)\le \id_\P\le \uex_{\id_\A }\left(\lex_{\mu \circ (\psi\times \id_\A)}(\varphi,\varphi)\right)\,.
\]
Therefore, applying $\uex_{\id_\A}$ to both sides of (\ref{eq:CIC2}) and using Proposition \ref{prop:inequality} yields
\[
\id_{\widehat \P}\le \uex_{\id_\A }\left(\lex_{\mu \circ (\psi\times \id_\A)}(\varphi,\varphi)\right) \le \uex_{\id_\A}(\lex_\psi(\varphi))\circ \varphi \le \uex_{\id_\A}(\id_{\widehat \P})=\id_{\widehat \P}\,,
\]
where the last equality follows from the fact that $\id_{\widehat \P}$ is an element of $\A$.
Hence $\varphi$ has a left inverse. A similar argument shows that it has right inverse and concludes the proof.

\begin{cor}\label{cor:CIC}
Let $\P$ be a Darboux complete ordered set and let $\A\subseteq \OP(\widehat\P)$ be a commutative completely integrally closed subgroup. Then $\BDar(\id_\A)$ is a commutative group. 
\end{cor}

\proof Let $\mu':\OP(\widehat\P)\times\OP(\widehat\P)\to\OP(\widehat\P)$ denote composition in reverse order i.e.\ $\mu'(\varphi,\varphi')=\varphi'\circ\varphi$, for all $\varphi,\varphi'\in \OP(\widehat\P)$. Since the restrictions of $\lex_{\id_\A}\circ \mu'$  and $\uex_{\id_\A}\circ \mu'$ to $\A\times \A$ coincides with $\mu\circ (\id_\A\times \id_\A)$, we obtain
\[
\lex_{\mu\circ (\id_\A \times \id_\A)} \le \lex_{\id_\A}\circ \mu'\le \uex_{\id_\A}\circ \mu'\le \uex_{\mu\circ (\id_\A \times \id_\A)}\,.
\]
Together with \eqref{eq:CIC1}, this implies the commutativity of the monoid $\Dar(\id_\A)$ which contains $\BDar(\id_\A)$.

\section{The Darboux Completion}\label{sec:completion}

\begin{rem}\label{rem:completion} Let $\O$ be an ordered set.
For each $x\in\O$, let $\delta_x:\O^\op\hto\widehat\emptyset$ be the partial function such that $\Dom(\delta_x)=\{x\}$ and $\delta_x(x)=+\infty$. Then $\lex_{\delta_x}(y)=+\infty$ if and only if $y\le x$. Let us denote $Y(x)=\lex_{\delta_x}$ for every $x\in \O$. If $f\in \O^\vee$ then $f(x)=+\infty$ if and only if $Y(x)\le f$.  Moreover, $f\le g$ in $\O^\vee$ if and only if $Y(x)\le f$ implies $Y(x)\le g$. In particular, $Y(x)\le Y(y)$ if and only if $x\le y$. Hence the assignment $x\mapsto Y(x)$ defines an order preserving embedding $Y:\O\to\O^\vee$ called the {\it Yoneda embedding of $\O$}.
\end{rem}

\begin{prop}\label{prop:completion}
Let $\O$ be an ordered set, let $\varphi:\O^\vee\hto \O^\vee$ be the identity function of the image  of the Yoneda embedding of $\O$ and let $\Dar(\O)$ denote the Darboux set of $\varphi$. Then
\begin{enumerate}[1)]
\item $\lex_\varphi = \id_{\O^\vee}$;
\item if $g\in\OP(\Dar(\O))$ restricts to the identity on $Y(\O)$, then $g=\id_{\Dar(\O)}$;
\item $\uex_\varphi(\O^\vee)\subseteq \Dar(\O)$;
\item the empty partial function of $\Dar(\O)$ is extremizable.
\end{enumerate}
\end{prop}

\proof Since $\id_{\O^\vee}$ restricts to $\varphi$ on $Y(\O)$, then $\lex_\varphi(f)\le f$ for every $f\in \O^\vee$. On the other hand, $Y(x)\le f$ implies $Y(x)=\lex_\varphi(Y(x))\le \lex_\varphi(f)$. Using Remark \ref{rem:completion}, this proves 1). 2) follows immediately from 1) and the definition of $\Dar(\O)$. Proposition \ref{prop:inequality} and 1) yield
\[
\uex_\varphi=\uex_\varphi\circ\lex_\varphi\le \uex_\varphi\circ\uex_\varphi\le \uex_{\varphi\circ\varphi}=\uex_\varphi
\]
which readily implies 3). Since $+\infty\le \uex_\varphi(+\infty)\le +\infty$, then $+\infty\in \Dar(\O)$. If $-\infty\neq \Dar(\O)$, then by Remark \ref{rem:completion} there exists $x\in \O$ such that $Y(x)\le \uex_\varphi(-\infty)$. By Lemma \ref{lem:complete} this implies that $x\le y$ for all $y\in \O$. Therefore, the empty partial function of $\Dar(\O)$ is extremizable, $\uex_\emptyset(\O)=+\infty$ and either $\lex_\emptyset(\O)$ is the function that takes value $-\infty$ on the complement of a set of cardinality at most one.

\begin{dfn}
Using the notation of Proposition \ref{prop:completion} and Remark \ref{rem:3.12}, we define the {\it Darboux completion of an ordered set $\O$} to be the ordered set
\[
\Dar'(\O)=\Dar(\O)\setminus \{\lex_\emptyset(\Dar(\O)),\uex_\emptyset(\Dar(\O))\}\,.
\]
\end{dfn}

\begin{cor}\label{cor:completion} The Darboux completion of an ordered set is Darboux complete.
\end{cor}

\proof Let $\O$ be an ordered set and let $\iota:\Dar(\O) \to \O^\vee$ be the inclusion. For any partial function $\psi:\Dar(\O)\hto\Dar(\O)$, Example \ref{ex:4.6} ensures that $\iota\circ \psi$ is extremizable. By Proposition \ref{prop:completion},
$\uex_\varphi\circ \uex_{\iota\circ\psi}$ and $\uex_\varphi\circ\lex_{\iota\circ\psi}$ are order preserving  functions in $\OP(\Dar(\O))$ that restrict to $\psi$ on $\Dom(\psi)$. On the other hand, $\lex_{\iota\circ\psi}\le\iota\circ  g\le \uex_{\iota\circ\psi}$ for any extension $g$ of $\psi$ to $\Dar(\O)$. Since $\uex_\varphi\circ \iota\circ g=g$, this implies that $\uex_\varphi\circ \lex_{\psi'} \le g\le \uex_\varphi\circ\uex_{\psi'}$ and thus $\psi$ is extremizable. This concludes the proof, since by construction $\Dar(\O)$ is canonically isomorphic to $\widehat{\Dar'(\O)}$.

\begin{rem}\label{rem:yoneda-convention}
From now on we use the Yoneda embedding to canonically identify $\O$ with a subset of $\Dar(\O)\subseteq \O^\vee$. In particular, this provides a canonical embedding of $\OP(\O)$ into the set of partial functions $\Dar(\O)\hto\Dar(\O)$.
\end{rem}

\begin{example}\label{ex:6.6}
We define the set of {\it real numbers} to be the Darboux completion $\R$ of the ordered set $\Q$ of rational numbers.  Moreover, $\Dar(\Q)$ is canonically identified with the set of {\it extended real numbers} $\widehat \R = \R\cup\{\pm\infty\}$. 
\end{example}

\begin{example}
Let $\Q_{>0}\subseteq \Q$ be the ordered set of positive rational numbers and let $\R_{>0}=\Dar'(\Q_{>0})$. Extending each function in  $\Dar'(\Q_{>0})$ by $+\infty$ to $\Q\setminus \Q_{>0}$ yields a canonical embedding of $\R_{>0}$ into $\R$ whose image consists of real numbers that are greater than $Y(0)$. Moreover, the composition of this embedding with the canonical embedding of $\Q_{>0}$ into $\Dar'(\Q_{>0})$ coincides with the restriction of the canonical embedding of $\Q$ into $\R$. Keeping in mind the above canonical identifications, it makes sense to write equalities such as $\Q_{>0}=\Q\cap \R_{>0}$. 
\end{example}

\begin{rem}\label{rem:restriction}
Let $\O$ be an ordered set, let $\P$ be a complete ordered set and let $\psi:\Dar(\O)\hto \widehat\P$ be an embedding with domain $Y(\O)$ and inverse $\psi':\widehat\P\hto \Dar(\O)$. Since $\uex_{\psi'}\circ\uex_\psi$ restricts to the identity on $Y(\O)$, then it equals $\id_{\Dar(\O)}$ by Proposition \ref{prop:completion}. Therefore, $\uex_\psi:\Dar(\O)\to \widehat \P$ is an embedding. In particular, it can attain the values $\pm\infty$ at most once which implies that $\uex_\psi$ restricts to an embedding $f:\Dar'(\O)\to \P$. By Remark \ref{rem:4.9} this implies that $\Dar'(\O)$ satisfies the same universal property of the Dedekind-MacNeille completion of $\O$ and is therefore canonically isomorphic to it. In particular, this shows that our definition of $\R$ is canonically isomorphic to the ordered set $\R'$ of Dedekind cuts of $\Q$. In fact, in this case it is easy to see directly that $\uex_{\psi'}:\widehat \R'\to \widehat \R$ is injective since it maps the cut associated to a rational number $x$ to $Y(x)$ and  $(\uex_{\psi'}(C))^{-1}(+\infty)=C$ for any irrational cut $C$. 
\end{rem}

\begin{prop}
There exists a canonical embedding $\alpha:\Aut(\O)\to\Aut(\Dar(\O))$. Moreover, $\alpha$ is a group homomorphism. 
\end{prop}

\proof Let $\varphi\in \Aut(\O)$. Using the convention of \ref{rem:yoneda-convention}, we may think of $\varphi$ as a partial function $\Dar(\O)\hto \Dar(\O)$. Then by Remark \ref{rem:restriction} 
\[
\lex_{\varphi^{-1}}\circ \uex_\varphi=\id_{\Dar(\O)}=\uex_\varphi\circ \lex_{\varphi^{-1}}\,.
\]
This implies that $\uex_\varphi$ is invertible and $\uex_\varphi\le \lex_\varphi \circ \lex_{\varphi^{-1}}\circ \uex_\varphi \le \lex_\varphi$.
Therefore, $\ex_\varphi\in \Aut(\Dar(\O))$. Let $\alpha(\varphi)=\ex_\varphi$ for all $\varphi\in \Aut(\O)$. Combining Remark \ref{rem:restriction} and Proposition \ref{prop:inequality}, we conclude that $\alpha$ is an injective group homomorphism and the Proposition is proved.

\begin{example}\label{ex:6.10}
Addition in $\Q$ defines an embedding $\lambda:\Q\to\Aut(\Q)$ such that $(\lambda(r))(s)=r+s$ for all $r,s\in \Q$. Composing with $\alpha$ we obtain an embedding $\beta:\Q\to\Aut(\widehat\R)$. Since every (order preserving) automorphism of $\widehat \R$ necessarily fixes $\pm \infty$, we have a canonical identification of  $\Aut(\widehat \R)$ with $\Aut(\R)$. In particular, $(\beta(x))(\pm\infty)=\pm \infty$ for all $x\in \Q$.
\end{example}

\begin{prop}\label{prop:6.11}
$\R$ is canonically isomorphic to $\BDar(\id_{\beta(\Q)})$.
\end{prop}

\proof Considering the embedding $\beta$ constructed in Example \ref{ex:6.10} as a partial function $\R\hto\OP(\widehat \R)$ (which is extremizable by Proposition \ref{prop:4.4}), we obtain order preserving functions $\lex_\beta,\uex_\beta:\R\to\OP(\widehat \R)$. The order preserving function $\ev_0:\BDar(\id_{\beta(\Q))})\to\R$ is surjective by Remark \ref{rem:restriction} since $\ev_0\circ \lex_\beta$ and $\ev_0\circ\uex_\beta$ both restrict to the identity on $\Q$. Since $\lex_\beta\circ \ev_0$ and $\uex_\beta\circ \ev_0$ both restrict to the identity on $\beta(\Q)$, they both equal to the identity on $\BDar(\id_{\beta(\Q)})$. Therefore $\ev_0$ is invertible with inverse $\ex_\beta$. 

\begin{rem} Combining Proposition \ref{prop:6.11} with
Corollary \ref{cor:CIC}, we conclude that $\R$ has a canonical structure of commutative group. Alternatively, this structure can be understood as follows. Let $+$ be the addition operation on $\Q$, thought of as a partial function $\R\times \R\hto \R$. Since $(\ex_\beta(r))(s)=r+s$ for all $r,s\in \Q$, we obtain
\begin{equation}\label{eq:addition}
\lex_+(x,y)\le (\ex_\beta(x))(y)\le \uex_+(x,y)
\end{equation}
for all $x,y\in \R$. On the other hand, for every $r\in \Q$ both $\ex_\beta(r)^{-1}\circ \uex_+(r,-)$ and $\ex_\beta(r)^{-1}\circ \lex_+(r,-)$ restrict to the identity of $\Q$. By Remark \ref{rem:restriction}, this implies $\lex_+(r,-)=\ex_\beta(r)=\uex_+(r,-)$ for all $r\in \Q$ and thus $\lex_\beta(x)\le\lex_+(x,-)\le \uex_+(x,-)\le \uex_\beta(x)$ for all $x\in \R$. Hence the inequalities of \eqref{eq:addition} are actually equalities for all $x,y\in \R$.
\end{rem}

\begin{rem}
A similar argument shows that the multiplication on $\Q_{>0}$ thought of as a partial function $\bullet: \R_{>0}\times\R_{>0}\hto\R_{>0}$ defines an partial function $\gamma:\R_{>0}\hto\Aut(\R_{>0})$ such that $\Dom(\gamma)=\Q_{>0}$ and $(\ex_\gamma(x))(y)=\ex_\bullet(x,y)$ for all $x,y\in \R_{>0}$. 
\end{rem}

\begin{thm}
$(\R_{>0},\ex_+,\ex_\bullet)$ is a semifield. 
\end{thm}

\proof Let $\psi:(\R_{>0})^3\hto\R_{>0}$ be the partial function with domain $(\Q_{>0})^3$ and such that $\psi(r,s,t)=r(s+t)$ for all $r,s,t\in \Q_{>0}$. Since $\ex_\bullet(r,\ex_+(s,t))=\psi(r,s,t)$ for all $r,s,t\in \Q_{>0}$, then
\begin{equation}\label{eq:distributivity-1}
\lex_\psi(x,y,z)\le \ex_\bullet(x,\ex_+(y,z)) \le \uex_\psi(x,y,z)
\end{equation}
for all $x,y,z\in \R_{>0}$. On the other hand, since $\gamma(s+t)$ agrees with both $\lex_\psi(-,s,t)$ and $\uex_\psi(-,s,t)$ on $\Q_{>0}$ for all $s,t\in \Q_{>0}$, they also agree on $\R_{>0}$. Using $(\gamma(s+t))(x)=\ex_\bullet(x,s+t)=(\ex_\gamma(x))(s+t)$, we obtain
\[
\lex_+(y,z)\le (\ex_\gamma(x))^{-1} \lex_\psi(x,y,z) \le (\ex_\gamma(x))^{-1} \uex_\psi(x,y,z) \le \uex_+(y,z)\,.
\]
Since $\ex_\bullet(x,\ex_+(y,z))=(\ex_\gamma(x))(\ex_+(y,z))$ for all $x,y,z\in \R_{>0}$, we conclude that the equalities of \eqref{eq:distributivity-1} are in fact equalities. It follows from the distributivity of $\bullet$ over $+$ on $\Q_{>0}$ that $\ex_+(\ex_\bullet(r,t),\ex_\bullet(r,t))=\psi(r,s,t)$ for all $r,s,t\in \Q_{>0}$. Hence, $\ex_\bullet$ distributes over $\ex_+$ on $\R_{>0}$ and the theorem is proved.

\begin{rem}
A standard argument shows that $\ex_\bullet$ can be canonically extended to an operation (which is not order-preserving) $\cdot$ on $\R$ in such a way that $(\R,\ex_{+},\cdot)$ is a field. With a slight abuse of notation, from now on we write $+$ for $\ex_{+}$. Since $(\beta(x))(\pm\infty)=\pm\infty$ for all $x\in \Q$, we set $x+(\pm \infty)=\pm\infty$ for all $x\in \R$. 
\end{rem}

\begin{rem}
For each ordered set $\O$, the set $\OP(\O,\R)$ inherits a canonical structure of $\R$-algebra with operations defined pointwise on $\O$. In particular, $f_1\le f_2$ implies $f_1+f_3\le f_2+f_3$ for any $f_1,f_2,f_3\in \OP(\O,\R)$ and $f_1f_3\le f_2f_3$ whenever $0\ge f_3$.
\end{rem}

\section{Limits and Integrals}\label{sec:limits}

\begin{dfn}
A {\it filter basis} on an ordered set $\O$ is a collection $F$ of nonempty subsets of $\O$ that is closed under finite intersections. To each filter basis $F$ of $\O$ we associate the partial function $\psi_F: \OP(\O,\R)\hto\widehat \R$ such that
\[
\Dom(\psi_F)=\bigcup_{\S\in F} \Dom(\psi_\S)\,,
\]
where $\psi_\S$ is defined as in Example \ref{ex:constant} and $\psi_F(f)=\psi_\S(f)$ for each $f\in \Dom(\psi_\S)$ and for each $\S\in F$.
\end{dfn}

\begin{dfn}\label{def:F-convergent}
Let $\O$ be a discrete set and let $F$ be a filter basis on $\O$. An order preserving function $f:\O\to\R$ is {\it $F$-convergent} if there exists $\lim_F(f)\in \widehat \R$ such that for every $\varepsilon>0$ there exists $\S\in F$ such that $f(x)\in [\lim_F(f)-\varepsilon,\lim_F(f)+\varepsilon]$ for all $x\in \S$.
\end{dfn}

\begin{thm}\label{thm:limit}
Let $\O$ be an discrete set and let $F$ be a filter basis on $\O$. An order preserving function $f:\O\to\R$ is $F$-convergent if and only if $f\in \Dar(\psi_F)$. Moreover, $\ex_{\psi_F}(f)=\lim_F(f)$ for all $f\in \Dar(\psi_F)$.
\end{thm}

\proof Assume that $f\in \Dar(\psi_F)$. Then by Remark \ref{rem:partition-1} for every $\varepsilon>0$ there exist $\S',\S''\in F$ such that $[\ex_{\psi_F}(f),\ex_{\psi_{S'}}(f)]\subseteq[\ex_{\psi_F}(f),\ex_{\psi_F}(f)+\varepsilon]$ and $[\ex_{\psi_{S''}}(f),\ex_{\psi_F}(f)]\subseteq [\ex_{\psi_F}(f)-\varepsilon,\ex_{\psi_F}(f)]$. Therefore, setting $\S=\S'\cap \S''$ we obtain
\[
\ex_{\psi_F}(f)-\varepsilon\le \lex_{\psi_\S}(f)\le f(x) \le \uex_{\psi_\S}(f)\le \ex_{\psi_F}(f)+\varepsilon
\]
for every $x\in \S$. Hence $f$ is $F$-convergent and $\lim_F(f)=\ex_{\psi_F}(f)$. Conversely, if $f$ is $F$-convergent, for every $\varepsilon>0$ there exists $\S\in F$ such that $\lim_F(f)-\varepsilon\le f(x)\le \lim_F(f)+\varepsilon$ for all $x\in \S$. Using Remark \ref{rem:constant}, this implies
\[
[\lex_{\psi_F}(f),\uex_{\psi_F}(f)]\subseteq [\lex_{\psi_\S}(f),\uex_{\psi_\S}(f)]\subseteq [\lim_F(f)-\varepsilon,\lim_F(f)+\varepsilon]
\]
for every $\varepsilon>0$. Hence, $f\in \Dar(\psi_F)$ and $\lim_F(f)=\ex_{\psi_F}(f)$.

\begin{rem}
In terms of the philosophy outlined in Section \ref{sec:introduction}, the functions in $\Dom(\psi_F)$, i.e. the functions that are constant on some element of $F$ are ``obviously $F$-convergent'' and $\psi_F$ is their ``obvious limit''. Feeding the machinery of Darboux calculus with this information results in a construction of general $F$-convergent functions that is alternative to the $(\varepsilon,\delta)$-definition given in \ref{def:F-convergent}.
\end{rem}

\begin{example}\label{ex:sequences}
Let $\N$ be the set of natural numbers with its usual order. Let $\O=|\N|$ and let $F=\{\N\setminus [1,n]\}_{n\in \N}$. Then $\OP(\O,\R)$ is the set of all sequences, $\Dom(\psi_F)$ is the set of sequences that are eventually constant and $\psi_F(f)$ is the function that to such a sequence assigns its obvious limit i.e.\ the only value that $f$ attains infinitely many times. Moreover, $\Dar(\psi_F)$ is precisely the set of all convergent sequences (including those converging to $\pm\infty$) and $\ex_{\psi_F}$ is their limit. 
\end{example}

\begin{example}\label{ex:continuous}
Let $\O=|\R|$, let $x_0\in \R$ and let $F$ be the collection of all subsets of the form $[x_0-\delta,x_0+\delta]$ for some $\delta>0$. Then $\OP(\O,\R)$ is the set of all real-valued functions of one real variable, $\Dom(\psi_F)$ is the subset of functions that are constant in a neighborhood of $x_0$ and $\psi_F(f)=f(x_0)$ for all $f\in \Dom(\psi_F)$. Moreover, $\Dar(\psi_F)$ is precisely the set of all functions that are continuous at $x_0$ and $\ex_{\psi_F}(f)=f(x_0)$ for all $f\in \Dom(\psi_F)$.
\end{example}

\begin{example}\label{ex:limits} In the notation of Example \ref{ex:continuous}, we could also consider $F$ to be the collection of all subsets of the form $[x_0-\delta,x_0+\delta]\setminus\{x_0\}$ for some $\delta>0$. Then $f\in \Dar(\psi_F)$ if an only if $f$ has limit at $x_0$ in which case $\ex_{\psi_F}(f)$ equals the limit. We leave the obvious variations leading to left and right limits to the reader.
\end{example}

\begin{dfn}
We denote by $\Int(\O)$ the ordered set of all intervals with endpoints in $\O$ with order given by $[a,b]\le [c,d]$ if and only if $c\le a\le b\le d$. We write $\int(\O)$ for the collection of nonempty subsets of $\O$ of the form $[x,z]\setminus\{x,z\}$, for some $[x,z]\in \Int(\O)$ (also ordered by inclusion). 
\end{dfn}

\begin{example}\label{ex:integrals} For any $J\in \Int(\R)$ let $m:\int(J)\to\R_{> 0}$ be the order preserving function defined by $m([x,z]\setminus \{x,z\})=z-x$ whenever $x\le z$ and $0$ otherwise. Given $J\in \Int(\R)$, let $\Par(J)$ be the collection of {\it partitions} of $J$ i.e.\ finite collections $P\subseteq \int(J)$ of mutually disjoint subsets such that $J\setminus \bigcup_{I\in P} I$ is finite. For each $I\in \int(J)$, let $\psi_I:\OP(|J|,\R)\hto\widehat \R$ be the partial function associated to the nonempty subset $I$ as in Example \ref{ex:constant}. In particular, the set of all bounded $\R$-valued functions on $J$ coincides with
\[
\O=B(\psi_J)=\bigcap_{I\in P} B(\psi_I)
\]
for any $P\in \Par(J)$. Let $\varphi_I:\B(\psi_I)\hto\R$ be the extremizable partial function defined as in Remark \ref{rem:constant} for each $I\in \int(J)$ and let $U_P,L_P:\O\to\R$ be the order preserving functions defined by
\[
U_P=\sum_{I\in P} m(I)\uex_{\varphi_I}\quad\textrm{ and }\quad L_P=\sum_{I\in P} m(I)\lex_{\varphi_I}\,.
\]
By Remark \ref{rem:constant} it is clear that for each $f\in \O$, $U_P(f)$ coincides with the usual {\it upper Darboux sums} of $f$ and $L_P(f)$ coincides with the usual {\it lower Darboux sums} of $f$ (see e.g.\cite{rudin}). On the other hand, for each $P\in \Par(J)$ consider the partial function $\varphi_P:\O\hto\R$ defined by
\begin{equation}\label{eq:integraldefinition}
\varphi_P(f)=\sum_{I\in P} m(I)\varphi_I(f)
\end{equation}
for each element $f$ of
\[
\Dom(\varphi_P)=\bigcap_{I\in P} \Dom(\varphi_I)\,.
\]
Since $\varphi_P$ is clearly encompassing, it is also extremizable by Corollary \ref{cor:encompassing}. Moreover,
\begin{equation}\label{integral}
\lex_{\varphi_P}\le L_P \le U_P \le \uex_{\varphi_P}
\end{equation}
as each term in the above chain of inequalities restricts to $\varphi_P$ on $\Dom(\varphi_P)$. Since each function in $\O$ attains only finitely many values, then $\Dom(\varphi_P)\cong \R^N$ for some integer $N$. In particular, Corollary \ref{cor:encompassing} ensures that if $\rho_P:\O\hto\Dom(\varphi_P)$ denotes the identity function on $\Dom(\varphi_P)$ then $\rho_P$ is extremizable. Therefore
\begin{equation}\label{integral-2}
\uex_{\varphi_P}\le \varphi_P\circ \uex_{\rho_P} = \sum_{I\in P} m(I)(\varphi_I\circ \uex_{\rho_P})\le \sum_{I\in P}m(I)\uex_{\varphi_I\circ \rho_P} = U_P\,.
\end{equation}
Combined with an analogous estimate for $\lex_{\varphi_P}$ and \eqref{integral}, \eqref{integral-2} shows that $\uex_{\varphi_P}=U_P$ and $\lex_{\varphi_P}=L_P$. Since $m(I)=m(I_1)+m(I_2)$ whenever $I\setminus (I_1\cup I_2)$ is finite and $I\subseteq I'$ implies $\varphi_I(f)=\varphi_{I'}(f)$ for each $f\in \Dom(\varphi_I')$, then for each $f\in \Dom(\varphi_P)\cap \Dom(\varphi_{P'})$
\[
\varphi_P(f)=\sum_{I\in P}m(I)\varphi_I(f)=\sum_{I\in P}\sum_{I'\in P'}m(I\cap I')\varphi_{I\cap I'}(f)=\sum_{I'\in\P'}m(I')\varphi_{I'}(f)=\varphi_{P'}(f)\,.
\]
Let $\varphi$ be the common extension of the compatible set $\{\varphi_P\}_{P\in \Par(J)}$. In particular, $\Dom(\varphi)$ is the set of {\it step functions on $J$} i.e.\ the set of all functions on $|J|$ that are constant on each interval of some partition of $J$. Combining \ref{rem:partition-2} with Corollary \ref{cor:encompassing} we conclude that $\varphi$ is encompassing and thus extremizable. By Remark \ref{rem:partition-3}, an order preserving function $g:\O\to\R$ restricts to $\varphi$ on $\Dom(\varphi)$ if and only if
\[
g\in\bigcap_{P\in\Par(J)}[\lex_{\varphi_P},\uex_{\varphi_P}]=\bigcap_{P\in\Par(J)}[L_P,U_P]\,.
\]
Hence, $\lex_\varphi$, $\uex_\varphi$ coincide with the lower and upper integral of $f$ on $J$, respectively. In particular, $\Dar(\varphi)$ coincides with the set of functions on $J$ that are integrable in the sense of Riemann and $\ex_\varphi$ is the Riemann integral. This is in fact the motivating example for the philosophy of Section \ref{sec:introduction}: to define the Riemann integral it is sufficient to feed the ``obvious'' definition for step functions (given by \eqref{eq:integraldefinition}) into the machinery of Darboux calculus to automatically obtain the correct general definition.
\end{example}

\begin{thm}[Linearity]\label{thm:additivity}
Let $\psi:\OP(\O,\R)\hto\R$ be a partial function such that
\begin{enumerate}[1)]
\item $\psi$ is encompassing;
\item $\Dom(\psi)$ is an $\R$-linear subspace of $\OP(\O,\R)$;
\item $\psi$ is an $\R$-linear transformation.
\end{enumerate}
Then for every $f_1,f_2\in \OP(\O,\R)$ and for every non-negative $a_1,a_2\in \R$
\begin{equation}\label{eq:linearity-1}
\uex_\psi(a_1f_1+a_2f_2)\le a_1\uex_\psi(f_1)+a_2\uex_\psi(f_2)
\end{equation}
and similarly
\begin{equation}\label{eq:linearity-2}
a_1\lex_\psi(f_1)+a_2\lex_\psi(f_2)\le \lex_\psi(a_1f_1+a_2f_2)\,.
\end{equation}
Moreover
\begin{equation}\label{eq:linearity-3}
-\lex_\psi(f)=\uex_\psi(-f)
\end{equation}
for every $f\in \OP(\O,\R)$. In particular, $\Dar(\psi)$ is an $\R$-linear subspace of $\OP(\O,\R)$ and $\ex_\psi$ is $\R$-linear.
\end{thm}

\proof  Since $\psi$ is encompassing, it is extremizable. By the additivity of $\psi$, the assignment $f_1\mapsto \uex_\psi(f_1+\eta_2)-\psi(\eta_2)$ coincides with $\psi$ on $\Dom(\psi)$ for each fixed $\eta_2\in \Dom(\psi)$. Therefore
\begin{equation}\label{eq:linearity-4}
\uex_\psi(f_1+\eta_2)\le \uex_\psi(f_1)+\psi(\eta_2)
\end{equation}
for every $f_1\in \OP(\O,\R)$ and for every $\eta_2\in \Dom(\psi)$. In particular
\[
\uex_\psi(f_1)+\psi(\eta_2)=\uex_\psi((f_1+\eta_2)+(-\eta_2))+\psi(\eta_2)\le \uex_\psi(f_1+\eta_2)
\]
and thus the inequality in \eqref{eq:linearity-4} is actually an equality. Therefore, for each $f_1\in \OP(\O,\R)$ the assignment $f_2\mapsto \uex_\psi(f_1+f_2)-\uex_\psi(f_1)$ coincides with $\psi$ on $\Dom(\psi)$. This proves \eqref{eq:linearity-1} when $a_1=a_2=1$. Since $\psi$ is compatible with scalar multiplication, $a$ is a positive real number and the assignment $f\mapsto a^{-1}\uex_\psi(af)$ coincides with  $\psi$ on $\Dom(\psi)$. Therefore, $\uex_\psi(af)\le a\uex_\psi(f)$ which in turn implies
\[
a\uex_\psi(f)=a\uex_\psi(a^{-1}(af))\le \uex_\psi(af)
\]
As a result, $a\uex_\psi(f)=\uex_\psi(af)$ for ever positive real number $a$ and for every $f\in \OP(\O,\R)$. This proves \eqref{eq:linearity-1} and \eqref{eq:linearity-2} is proved similarly. Since $\psi$ is odd, the assignments  $f\mapsto -\lex_\psi(-f)$ and $f\mapsto -\uex_\psi(-f)$ both restrict to $\psi$ on $\Dom(\psi)$ and thus
\begin{equation}\label{eq:14}
\lex_\psi(f)\le -\uex_\psi(-f)\le -\lex_\psi(-f)\le \uex_\psi(f)
\end{equation}
for every $f\in \OP(\O,\R)$. The first inequality in \eqref{eq:14} implies $\uex_\psi(-f)\le -\lex_\psi(f)$ while the last inequality of \eqref{eq:14} implies $-\lex_\psi(f)\le \uex_\psi(-f)$ and thus \eqref{eq:linearity-3}. The last statement is a straightforward consequence of \eqref{eq:linearity-1}-\eqref{eq:linearity-3}.

\begin{example}
Specializing Theorem \ref{thm:additivity} to Examples \ref{ex:sequences}-\ref{ex:integrals} we obtain the well-known linearity theorems for limits of sequences, continuous functions, limits of functions of real variable and integrals. 
\end{example}


\begin{bibdiv}
\begin{biblist}

\bib{edalat}{article}{
   author={Edalat, Abbas},
   author={Lieutier, Andr{\'e}},
   title={Domain theory and differential calculus (functions of one
   variable)},
   journal={Math. Structures Comput. Sci.},
   volume={14},
   date={2004},
   number={6},
   pages={771--802},
}

\bib{fuchs}{book}{
   author={Fuchs, L.},
   title={Partially ordered algebraic systems},
   publisher={Pergamon Press, Oxford-London-New York-Paris; Addison-Wesley
   Publishing Co., Inc., Reading, Mass.-Palo Alto, Calif.-London},
   date={1963},
}

\bib{maclane}{book}{
   author={Mac Lane, Saunders},
   title={Categories for the working mathematician},
   series={Graduate Texts in Mathematics},
   volume={5},
   edition={2},
   publisher={Springer-Verlag, New York},
   date={1998},
}


\bib{rudin}{book}{
   author={Rudin, Walter},
   title={Principles of mathematical analysis},
   edition={3},
   note={International Series in Pure and Applied Mathematics},
   publisher={McGraw-Hill Book Co., New York-Auckland-D\"usseldorf},
   date={1976},
}


\bib{taylor}{book}{
   author={Taylor, Paul},
   title={Practical foundations of mathematics},
   series={Cambridge Studies in Advanced Mathematics},
   volume={59},
   publisher={Cambridge University Press, Cambridge},
   date={1999},
}

\bib{taylor-lambda}{article}{
   author={Taylor, Paul},
   title={A lambda calculus for real analysis},
   journal={J. Log. Anal.},
   volume={2},
   date={2010},
   pages={Paper 5, 115},
   issn={1759-9008},
}

\bib{univalent}{book}{
   author={The Univalent Foundations Program},
   title={Homotopy type theory---univalent foundations of mathematics},
   publisher={The Univalent Foundations Program, Princeton, NJ; Institute
   for Advanced Study (IAS), Princeton, NJ},
   date={2013},
}



\end{biblist}
\end{bibdiv}
\end{document}